\documentclass[a4paper,12pt]{article}
\usepackage{amsmath,amssymb,euscript,amsthm}
\newcommand{\1}{1\!\!\!\;{\rm I}}
\newcommand{\mbR}{{\mathbb R}}
\newcommand{\wt}{\widetilde}

\newcommand{\ov}{\overline}
\newcommand{\loc}{\mathop{\rm loc}}
\newcommand{\vf}{\varphi}
\newcommand{\ve}{\varepsilon}
\newcommand{\cF}{\mathcal F}
\newcommand{\rank}{\mathop{\rm rank}}
\newcommand{\pt}{\partial}
\newcommand{\otau}{\mathop{\tau}\limits^\circ}
\newcommand{\Ps}{\mathsf{P}}

\begin{document}\large

 \centerline{ \textsc{\Large Stochastic flows with reflection}}
 \centerline{\textsf{ Andrey PILIPENKO}}

{\small Institute of Mathematics of Ukrainian Academy of Sciences,
Dept. of Stochastic Processes,  3 Tereschenkovskaya Str., 01601
Kiev-4, UKRAINE. ( apilip@imath.kiev.ua)}

\vskip 5pt

{\small \textbf{Summary:} Some topological properties of
stochastic flow $\varphi_t(x)$ generated by stochastic
differential equation in a ${\mathbb R}^d_+$ with normal
reflection at the boundary are investigated. Sobolev
differentiability in initial condition is received. The absolute
continuity of the measure-valued process $\mu\circ\varphi_t^{-1}$,
where $\mu\ll\lambda^d,$ is studied. }

 \vskip 5pt

Flows generated by SDEs in Euclidean space is a well-studied topic
nowadays. It is well known for example (cf. \cite{Ku} and ref.
therein) that if the coefficients of SDE  are Liphitzian then the
SDE generates a flow of homeomorphisms, if coefficients are of the
class $C^{n+\ve}$ then SDE generates $C^n$-flow of
diffeomorphisms, equations for derivatives are obtained by formal
differentiation of the SDE etc.

Note that the similar questions for SDEs with reflection is much
harder to answer. Even the problems about coalescence of two
reflecting Brownian motions \cite{CL, CL1, Sheu, BBC}  or
differentiability of the Brownian reflecting flow ($\sigma(x)
=const$) \cite{Burdzy_diff, A} need accurate and non-trivial
considerations.

The article below was  published in Reports of Ukrainian Nat.Acad.
of Sci. \cite{P_DAN} (2005). Only some new references or minor
remarks are added.

Assume that functions $a_k : \mbR^d_+ \to \mbR^d$ satisfy the
Lipschitz condition. Here $\mbR^d_+=\mbR^{d-1}\times[0,\infty)$.
Consider an SDE  in $\mbR^d_+$ with normal reflection from the
boundary:
$$
\begin{cases}
d\vf_t(x)=a_0(\vf_t(x))dt+\sum^m_{k=1}a_k(\vf_t(x))dw_k(t)+\\
\ \ \ \ \ \ \ \ +\ov{n}\xi(dt,x), \ t\in[0,T],\\
\vf_0(x)=x, \ \xi(0,x)=0, \ x\in\mbR^d_+,
\end{cases}
\eqno(1)
$$
where $ \{w_k(t), k=1,\ldots,m\}$ are independent Wiener
processes, $\ov{n}=(0,\ldots,0,1)$ is a normal to hyperplane
$\mbR^{d-1}\times\{0\},$  for each fixed $x\in\mbR^d_+$ a process
$\xi(t,x)$ is non-decreasing in $t$, and
$$
\xi(t,x)=\int^t_0\1_{\{\vf_s(x)\in\mbR^{d-1}\times\{0\}\}}\xi(ds,x),
$$
i.e. $\xi(t,x)$ is increasing only on those instants of time when
$\vf_t(x)\in\mbR^{d-1}\times\{0\}.$ Lipschitz property of the
coefficients ensures the existence and the uniqueness to the
solution of  (1), cf. \cite{Tanaka}.

 \vskip5pt
  {  \textbf{1. Existence of continuous modification.}}

 {\bf Theorem   1 \cite{P1}}. {\it
There exists  a modification of the processes
$\varphi_{t}(x),\xi(t,x)$ (it will be denoted in the same way)
such that

1) for any  $x\in\mbR^d_+$, the pair $(\varphi_{t}(x),\xi(t,x)), \
t\geq 0$, is a solution of  (1);

2) for any $\omega\in\Omega$ processes $\varphi_{t}(x),\xi(t,x)$
are continuous in a pair of arguments $(t,x), t\geq 0,
x\in\mbR^d_+.$ }

 The Theorem   1 is proved in a way similar to the corresponding proof
 used for the solution of  SDE without reflection, cf. \cite{Ku},
 with the use of  Kolmogorov's theorem on existence of continuous modification.

It will be assumed further that $\varphi_{t}(x),\xi(t,x) $ are
already continuous.
 \vskip5pt
   \textbf{2. The joint motion of solutions started from different initial points.} It is well known
\cite{Ku} that a solution of an SDE (without reflection) generates
a flow of diffeomorphisms. However, the injectivity for reflecting
flow can be failed as the following example shows.

{\bf Example 1.} Let $d=1, m=1, a_0=0, a_1=1,$ i.e. $\varphi_t(x)$
is the reflected Brownian motion in $\mbR^1_+$  started from
$x\geq 0:$
$$
\varphi_t(x)=x+w(t)+\xi(t,x),x\geq0.
$$
 It is easy to see that $\varphi_t(x), \xi(t,x)$ is of the form
$$
\varphi_t(x)=\begin{cases}
w(t)-\min_{0\leq s\leq t}w(s), \ x=0,\\
w(t)+x, \ x>0 \ \mbox{и} \ \tau(x)\geq t,\\
\varphi_t(0), \ x>0 \ \mbox{и} \ \tau(x)<t,
\end{cases}
$$
$$
\xi(t,x)=\begin{cases}
-\min_{\tau(x)\leq s\leq t}w(s), \tau(x)<t,\\
0, \ \tau(x)\geq t,
\end{cases}
$$
where $\tau(x)$ is a moment, when the process $x+w(t)$ gets zero
for the first time.

In other words, $\varphi_t(x)$ is moving as $x+w(t)$ before
hitting  0, and then a motion of $\varphi_t(x)$ coincides with the
reflected Brownian motion
   $\varphi_t(0)$ started from zero.

The similar situation takes place in multi-dimensional space.

{\bf Theorem 2. \cite{P1}}  {\it Denote by
$\tau(x)=\inf\left\{t\geq 0 : \
\vf_t(x)\in\mbR^{d-1}\times\{0\}\right\}$ the moment of the first
hitting the hyperplane $\mbR^{d-1}\times\{0\}$ by a solution
started from $x\in\mbR^{d}_+.$

Then there exists a set $\Omega_0$ of probability 1 such that for
all $\omega\in\Omega_0$ the following statements hold true:

1) for all  $x,y\in\mbR^d_+, x\ne y$ and
$t<\max\{\tau(x),\tau(y)\}$ the inequality  $
\varphi_t(x)\ne\varphi_t(y) $ is satisfied;

2) for any  $x\in\mbR^d_+$  there exists
$y=y(x,\omega)\in\mbR^{d-1}\times\{0\},$  such that
 $ \varphi_{\tau(x)}(x)=\varphi_{\tau(x)}(y) $ if
 $\tau(x)<\infty.$ Moreover,
$$\varphi_t(x)=\varphi_t(y) \  \mbox{for}\  t\geq\tau(x).
$$ }
{\bf Remark.}  Informally this theorem can be formulated in the
following way. A particle started from a point
$x\in\mbR^{d-1}\times(0,\infty)$ does not hit any other particle
before getting the hyperplane $\mbR^{d-1}\times\{0\}.$  At the
instant  $\tau(x)$ it coalesces with some other particle, which
started from $\mbR^{d-1}\times\{0\}.$ After this both particles
moves together.

\vskip5pt
  \textbf{3. Characterization of inner and boundary points of random set $\vf_t(\mbR^d_+).$}

{\bf Theorem 3.\cite{P_Hausd}}
{\it  For almost all  $\omega$ and all $t\in[0,T]$ the following
equality of random sets takes place
$$
\partial\vf_t(\mbR^d_+) =\vf_t(\partial\mbR^d_+) = \vf_t\{x\in\mbR^{d}_+: \tau(x)\leq t\},
$$
where  $\tau(x)=\inf\{s\geq0: \vf_s(x)\in\mbR^{d-1}\times\{0\}\}$
is the moment of the first hitting the hyperplane
$\mbR^{d-1}\times\{0\}$ by the solution started from $x$.

Moreover, for all  $R>0$  Hausdorff measure $H^{d-1}$ of the set
$\partial\vf_t(\mbR^d_+)\cap \{x\in\mbR^d_+ \ : \ \|x\|\leq R\}$
is finite. }

\vskip5pt  \textbf{4. Differentiability with respect to initial
condition.}

As in  Example 1, there is no reasons to expect that a solution of
(1) is  continuously differentiable in   $x$ even when
coefficients of the SDE are infinite differentiable. However, it
can be proved that for any $t$ the mapping $x\to \vf_t(x)$ belongs
to a Sobolev space
$\mathop{\cap}\limits_{p>1}W^1_{p,\loc}(\mbR^{d}_+,\mbR^d).$

The equations for  $\nabla\vf_t$ are not classical equations of
stochastic analysis. We need the following definition.

{\bf Definition  1.} Let $w_1(t),\ldots,w_m(t)$ be independent
Wiener processes,
$\cF_t=\sigma(w_k(s), k=\ov{1,m}, s\leq t),$  $ a_k:
\mbR^l\times\mbR^p\to\mbR^l,  \ b_k: \mbR^l\times\mbR^p\to\mbR^p,
\ k=0,\dots,m, $ and $x_t$ be a continuous  $\cF_t$-measurable
stochastic process. Consider a random measure-valued process
$\nu(t)=\delta_0{1\!\!\!\;{\rm I}}_{\{x(t)=0\}},$ where $\delta_0$
is a probability measure on $\mbR$, assigned unit mass to a point
 zero.

A pair $(y_t,z_t)$ of $\cF_t$-adapted processes satisfies the
equation
$$
\begin{cases}
dy_t=a_0(y_t,z_t)dt+\sum^m_{k=1}a_k(y_t,z_t)dw_k(t)-y_{t-}d\nu(t),\\
dz_t=b_0(y_t,z_t)dt+\sum^m_{k=1}b_k(y_t,z_t)dw_k(t), t\geq0,
\end{cases}
\eqno(2)
$$
if:

1) $y_t, t\geq0$ has cadlag trajectories;

2) $z_t,t\geq0$ has continuous trajectories;

3)
$z_t=z_0+\int^t_0b_0(y_s,z_s)ds+\sum^m_{k=1}\int^t_0b_k(y_s,z_s)dw_k(s),
t\geq0$ a.s.;

4) for almost all $\omega$ the set $\{t\geq0: x_t=0\}$ is
contained in $\{t\geq0: y_t=0\};$

5) for any stopping time  $\tau$ such that $x_{\tau}\ne 0$ a.s.,
the following equality holds true
$$
y_t=y_{\tau}+\int^t_{\tau}a_0(y_s,z_s)ds+\sum^m_{k=1}\int^t_{\tau}
a_k(y_s,z_s)dw_k(s)
$$
for all  $t\in[\tau,{\otau}),$  ${\otau}=\inf\{t\geq\tau:
x_t=0\}.$

{\bf Theorem 4. \cite{P_Diff} } {\it Assume that  functions  $a_k,
b_k, k=\ov{0,m}$ satisfy Lipschitz condition. Then there exists a
unique solution of (2) for any non-random initial condition
$(y_0,z_0)$. }

{\bf Theorem 5. } {\it I. \textbf{\cite{P_Hausd}} If functions
$a_k: \mbR^d_+\to\mbR^d, k=\ov{0,m}$ satisfy Lipschitz condition
then for a.a. $\omega$ a mapping  $
\mbR^d_+\ni{u}\to\varphi_t(u)\in\mbR^d $ belongs to the space
$\mathop{\cap}\limits_{p>1}W^1_{p,\loc}(\mbR^{d}_+,\mbR^d)$ for
a.a. $t\geq 0.$

II. \textbf{\cite{P_Diff}} Assume that functions  $a_k:
\mbR^d_+\to\mbR^d, k=\ov{0,m}$, are continuously differentiable
and their derivatives are bounded. Suppose also that $
\sum^m_{k=1}(a_{k,d}(x))^2>0$ for all $x\in\mbR^{d-1}\times\{0\}$,
where $a_{k,d}$ is the $d$-th coordinate of a function
$a_k=(a_{k,1}, \dots, a_{k,d})^T.$

Then the Sobolev derivative $\nabla\vf_t(x)$ satisfies the SDE
$$
\begin{cases}
d\nabla\vf_t(x)=\nabla
a_0(\vf_t(x))\nabla\vf_t(x)dt+\sum\limits^m_{k=1}\nabla
a_k(\vf_t(x))\nabla\vf_t(x)
dw_k(t)-\\ \ \ \ \ \ \ \ \  \ \ \ \ \ \ -P\nabla\vf_{t-}(x)n(dt,x),\\
\nabla\vf_0(x)=\1,
\end{cases}
\eqno(3)
$$
where  $\1$ is an identity matrix,  $P$ is a matrix corresponding
to the orthoprojection on the $d$-th coordinate of the space
$\mbR^d,$ $n(dt,x)$ is a point random measure such that
$n(\{t\},x)=1$ iff $\vf_t(x)$ belongs to the hyperplane
$\mbR^{d-1}\times\{0\}.$ }

{\bf Remark.} Equation  (3) is understood in the sense of
Definition   1. In this case we take the $d$-th coordinate of
$\vf_t(x)$ as   $x_t$, the $d$-th row of $\nabla\vf_t(x)$ as
$y_t$, the first $(d-1)$ rows of $\nabla\vf_t(x)$ and $\vf_t(x)$
as the process  $z_t$.

{\bf Remark.} The process  $\nabla\vf_t(x)$ can be chosen
measurable in $t,x,\omega.$

{\bf Remark.} The similar result   for constant diffusion
coefficient was obtained in \cite{A}. Moreover, it was proved that
for all $x$ and a.a. $\omega$ the mapping $\vf_t$ is  continuously
differentiable in some neighborhood of $x.$

Let us compare Sobolev differentiability and usual
differentiability of the mapping $\vf_t(\cdot,\omega).$ It is well
known that if the diffusion matrix is a constant then for a.a.
$\omega$ and all $t$ the mapping $x\rightarrow\vf_t(x,\omega)$
satisfies Lipschitz condition. Therefore $\vf_t(x)$ is
differentiable for $\lambda^d$-a.a. $x\in G$ by Rademacher's
theorem \cite{Fed}. Since the usual and Sobolev derivatives are
equal (if they exist), so equation for usual derivative coincides
with that for Sobolev. It is not difficult to prove that the usual
derivatives exist not only for a.a. $x$, a.a. $\omega,$ but for
all $x$ and a.a. $\omega$. However
 almost everywhere local continuous differentiability is not evident.

It should be noted that \cite{A} does not imply that for a.a.
$\omega$ the mapping $x\rightarrow\vf_t(x,\omega)$  is
continuously differentiable. Really, for the process from Example
1:

$$\Ps( x\rightarrow\vf_t(x)\ \mbox{is continuously
 differentiable})=0$$

 but for each $ x_0>0:$

\noindent $  \Ps( x\rightarrow\vf_t(x)$ is continuously
 differentiable in some neighborhood of $x_0)=1!$

The fact that  $\mbR^d_+\ni x\rightarrow\vf_t(x)$ is not
continuously
 differentiable seems to be typical, because $\rank\nabla\vf_t(x)\leq d-1$ if $\tau(x)\leq
 t$ and $\vf_t(x),\ \tau(x)> t$, coincides with the flow without
 reflection, so $\det\nabla\vf_t(x),\ \tau(x)> t, \|x\|\leq r, $ is separated from
 zero for any $r>0.$

\vskip5pt
 \textbf{5. Absolute continuity of image-measures
  driven by $\vf_t.$} Let  $\mu$ be a finite measure in
$\mbR^d_+$ which is absolute continuous w.r.t. Lebesgue measure.
Consider a measure-valued process  $\mu_t=\mu\circ \vf_t^{-1}, \
t\geq 0.$

Let us introduce a random set $O_t(\omega)=\{x\in\mbR^d_+:
t<\tau(x)\},$ where $\tau(x)=\inf\{s\geq0:
\vf_s(x)\in\mbR^{d-1}\times\{0\}\}$ is the moment of the first
hitting the hyperplane $\mbR^{d-1}\times\{0\}$ by the process
$\vf_s(x)$.

{\bf Theorem 6. \cite{P_Trans}} {\it For a.a. $\omega$ and every
$t\geq0$ a measure $\mu_t$ is represented as a sum of orthogonal
measures $ \mu_t= \mu\big|_{O_t}\circ\vf^{-1}_t+
\mu\big|_{\mbR^d_+\setminus O_t}\circ\vf^{-1}_t,$  such that

a) the first  measure is absolute continuous w.r.t.
$d$-dimensional Lebesgue measure and the second one is singular;

b) the support of measure $\mu\big|_{\mbR^d_+\setminus
O_t}\circ\vf^{-1}_t$ is contained in the set
$\vf_t(\mbR^{d-1}\times\{0\})$ of the  $\sigma$-finite
$(d-1)$-dimensional Hausdorff measure $H^{d-1}.$ }

The proof of the first part of the theorem follows from \cite{BH},
and the second part follows from Theorem 3.

The next theorem gives a sufficient condition that ensures the
absolute continuity of $ \mu\big|_{\mbR^d_+\setminus
O_t}\circ\vf^{-1}_t $ with respect to  $
H^{d-1}\big|_{\pt\vf_t(\mbR ^d_+)}, $ which is  the restriction of
$H^{d-1}$ to the set ${\pt\vf_t(\mbR ^d_+)}$.

{\bf Theorem 7. \cite{P_Trans}} {\it Assume that for $\mu$-a.a.
$x\in\mbR^d_+:$
$$
\Ps(\rank \nabla\vf_t(x)\geq d-1, t\geq 0)=1. \eqno(4)
$$
Then with probability 1 the absolute continuity
$$
\mu\big|_{\mbR^d_+\setminus O_t}\circ\vf^{-1}_t\ll H^{d-1}\big|_{
\pt\vf_t(\mbR^d_+)}  \eqno(5)
$$
holds for all  $t\geq 0.$

Here $\mu\big|_{\mbR^d_+\setminus O_t}, \ H^{d-1}\big|_{
\pt\vf_t(\mbR^d_+)}$ are restrictions of measures  $\mu, \
H^{d-1}$ to the sets  ${\mbR^d_+\setminus O_t}, \
\pt\vf_t(\mbR^d_+)$, respectively.
 }

 {\bf Remark.} In contrast to Hausdorff measure $H^{d-1}$ in
$\mbR^d,$ its restriction to the set
$\pt\vf_t(\mbR^d_+)=\vf_t(\mbR^d_+\setminus O_t)$ is a
$\sigma$-finite measure  (Theorem   3). So the notion of absolute
continuity in (5) does not require any specification.

The proof is provided by using the co-area formula \cite{Fed}
similarly to the case   $m=n$, cf. \cite{BH, Bog}

The verification of Theorem 7 conditions is quite difficult. Let
us give more simple sufficient conditions that ensure (4). Assume
that functions $a_k, 0\leq k\leq m$, are continuously
differentiable.

Denote by  $U_{st}(x), s\leq t$,  a solution of the following
linear SDE:
$$
\begin{cases}
dU_{st}(x)=\nabla a_0(\vf_t(x))U_{st}(x)dt+\sum^m_{k=1}\nabla
a_k(\vf_t(x))U_{st}(x)dw_k(t),
t\geq s, \\
U_{ss}(x)=\1.
\end{cases}
$$
Let us represent a random set $A=A_t(x)=\{s\in[0,t]:
\vf^d_s(x)>0\}$ as a disjoint union of random intervals $
A=[\alpha_0(x),\beta_0(x))\cup(\alpha_1(x),
\beta_1(x)]\cup\mathop{\bigcup}\limits^\infty_{k=2}(\alpha_k(x),\beta_k(x)),
$ where $\alpha_0(x)=0, \beta_1(x)=t.$  Let $P$ be the same as in
Theorem  5.

{\bf Theorem 8. \cite{P_Trans}} {\it  Let the conditions of the
second part of Theorem 5 be satisfied, $t\geq 0.$  Assume that for
all $k\geq 0$ and a.a. $x\in U$
$$
\Ps\left(\rank((1-P)U_{\alpha_k\beta_k}(x)(1-P))=d-1\right)=1.
\eqno(6)
$$
Then  relation  (4) holds true. }

{\bf Corollary.} Assume that for all   $x\in \mbR^d_+, s\leq t:$
\newline $\Ps\left(\vf_t^d(x)=0,  \det\|\wt U_{st}(x)\|=0\right)=0, $
where $\wt U_{st}(x)$ is the  matrix getting out from $U_{st}(x)$
by deleting last raw and last column. Then (6) fulfills.

Observe that the assumption  of the Corollary is the requirement  of
non-hitting zero by two-dimensional Ito process, and this is usually
easier to check than (6).

 {\bf Example 2.} Let $d=2$, then assumptions of the Corollary to the
 Theorem 8 are satisfied if for  all  $x\in \mbR^2_+, y\in  \mbR^2, y\ne 0$
 vectors
 $(a_{k,1}(x))_{1\leq k\leq m}$ and  $(\nabla_ya_{k,2}(x))_{1\leq k\leq
 m}$, are linear independent.

\vskip5pt
 \textbf{6. Flows generated by SDE with reflection in arbitrary set.} Let
   $G\subset \mbR^d$ be a closed set with smooth boundary such that the following
    SDE with normal reflection on the boundary of  $G$ has a unique
    strong solution defined for all $x\in G,
t\geq 0:$
$$
\begin{cases}
d\vf_t(x)=a_0(\vf_t(x))dt+\sum^m_{k=1}a_k(\vf_t(x))dw_k(t)+\\
\ \ \ \ \ \ \ \ +\ov{n}(\vf_t(x))\xi(dt,x), \ t\geq 0, \ x\in G,\\
\vf_0(x)=x, \ \xi(0,x)=0, \
\xi(t,x)=\int^t_0\1_{\{\vf_s(x)\in\partial G\}} \xi(ds,x),
\end{cases}
$$
where  $\ov{n}(x)$ is the inward normal at a boundary point
$x\in\partial G$, $\xi(t,x)$ is continuous and non-decreasing in
$t$ process for every fixed $x\in G.$

Assume that for any  $x\in\pt G$ there exists a neighborhood
  $O(x)$ and $C^2$-diffeomorphism  $\alpha_x$
which transforms  the set $O(x)\cap  G$ into $\{x\in \mbR^d \ : \
\|x\|\leq 1, x^d\geq 0\}$ in such a manner that
$\nabla\alpha_x(y)={\ov n}(y)=(0,\dots,0,1), \ y \in \pt G\cap
O(x).$

Applying a localization of solutions,    all
 statements on flows in half-space can be easily generalized (of course
 with natural changes) to the case of SDE in $G.$ For example,  relation (6) will be of the form:
 $$
 \Ps\left(\rank(P(\vf_{\beta_k})U_{\alpha_k\beta_k}(x)P(\vf_{\alpha_k}))=d-1\right)=1, \eqno(7)$$
where $P(x)$ is orthoprojection on the orthogonal complement to
$\ov{n}(x)$ where $x\in \pt G.$

{\bf Example 3.} Let $G=\{x\in \mbR^2 : \|x\|\leq 1\}$ be a unit
disk, $\vf_t(x), t\geq 0, x\in G$ be a Brownian motion in  $G$
with reflection on the boundary.
 I.e. $\vf_t(x)$ is a solution of the SDE
$$
\begin{cases}
d\vf_t(x)=dw(t)+\ov{n}(\vf_t(x))\xi(dt,x), \ t\geq 0,\\
\vf_0(x)=x, \ \xi(0,x)=0, \ x\in G,
\end{cases}
\eqno(8)
$$
where  $w(t)$ is a two-dimensional Wiener process.

Let us describe inner and boundary points of the set $\varphi_t(G)$.
Now, the stopping time $\tau(x)$ from Theorem 2 is of the form
$\tau(x)=\inf\{t\geq 0 \ : \ x+w(t)\in \pt G\}$. So, the set of inner points of $\varphi_t(G)$ is equal to
 $\{ x + w(t) :
x\in G, t<\tau(x)\}.$

Introduce a stopping time  $\sigma=\inf\{t\geq 0 \ : \
\|w(t)\|=2\}$. Observe that  $\sup_x\tau(x)\leq \sigma$. Therefore
the analogues of Theorems 2,3 imply that for every $t\geq \sigma$
the random set  $\varphi_t(G)$ coincides with the nowhere dense
set $\varphi_t(\pt G)$ of finite Hausdorff measure  $H^1$.
Moreover,
$$
 \ \Ps\left( \forall t\geq \sigma \ \forall x,
\|x\|<1 \ \exists y\in \pt G, \ x\ne y :
\vf_t(x)=\vf_t(y)\right)=1.
$$
It is interesting to compare this result with \cite{CL}, where it
is proved that any two solutions of (8)  started from different
initial points of  $G$ never meet each other with probability 1,
that is
$$
  \forall x,y\in G, \ x\ne y : \ \ \Ps\left(\exists t\geq 0 \ : \ \vf_t(x)=\vf_t(y)\right)=0.
$$

Note that in this example  $U_{st}(x)$ is an identity matrix, so
condition (7) is obviously satisfied. Thus for any absolute
continuous measure $\mu$ on $G$ we have the absolute continuity
 $$\mu\big|_{G\setminus
O_t}\circ\vf^{-1}_t\ll H^{1}\big|_{ \pt\vf_t(G)} \eqno(9)$$ with
probability 1. In particular,  for  a.a. $\omega$ and all $t\geq
\sigma$:
$$
\mu\circ\vf^{-1}_t\ll H^{1}\big|_{ \pt\vf_t(G)}.
$$

Observe that if $G$ is not a unit disk but any domain with "nice"\
boundary, and the boundary does not contain any two perpendicular
segments, then (7) and so (4),(9) are also satisfied. Exactly the
same condition on the boundary appears in \cite{Burdzy_diff}.
Moreover, it can be easily  shown that if the boundary contains
two perpendicular segments then (4),(7), and (9) are false.

\end{document}